\documentclass[12pt]{article}
\usepackage{amsmath,amsthm,amsfonts,amssymb}
\begin{document}
\title{Topological Quantum Field Theories and\\Operator Algebras}
\author{{\sc Yasuyuki Kawahigashi}
\footnote{The author was supported in part by JSPS Grants.}\\
Department of Mathematical Sciences\\
University of Tokyo, Komaba, Tokyo, 153-8914, Japan\\
e-mail: {\tt yasuyuki@ms.u-tokyo.ac.jp}}
\date{}
\maketitle

\section{Introduction}
\label{intro}

We have seen much fruitful interactions between 3-dimensional
topology and operator algebras since the stunning discovery 
of the Jones polynomial for links
\cite{J2} arising from his theory of subfactors \cite{J1}
in theory of operator algebras.  In this paper, we review
the current status of theory of ``quantum'' topological invariants of
$3$-manifolds arising from operator algebras.  
The original discovery of topological invariants arising from
operator algebras was for knots and links, as above, rather
than $3$-manifolds, but here we concentrate on invariants for
$3$-manifolds.  On the way of studying such topological
invariants, we naturally go through topological invariants
of knots and links.  From operator algebraic data, we construct
not only topological invariants of $3$-manifolds, but also
topological quantum field theories of dimension 3, in the 
sense of Atiyah \cite{A}, as the title of this paper shows, but
for simplicity of expositions, we consider mainly complex
number-valued topological invariants of oriented compact
manifolds of dimension 3 without boundary.

All the constructions of such topological invariants we discuss 
here are given in the following steps.
\begin{enumerate}
\item Obtain combinatorial data arising from representation
theory of an operator algebraic system.
\item Realize a manifold concretely using basic building blocks.
\item Multiply or add the complex numbers appearing in the data
in Step 1, in a way specified by how the basic building blocks
are composed in Step 2, and compute the resulting complex number.
\item Prove that the complex number in Step 3 is independent of
how the basic building blocks are composed, as long as the
homeomorphism class of the resulting manifold is fixed.
\end{enumerate}

In Step 1, the prototype of the representation theory for operator
algebras is 
the one for finite groups.  That is, for a finite group $G$,
we consider representatives of unitary equivalences classes
of irreducible unitary representations.  This finite set has
an algebraic structure arising from the tensor product
operation of representations, and it produces combinatorial
data such as fusion rules and $6j$-symbols.
In our setting, we work on some form of representation
theory of operator algebraic systems analogous to this classical
representation theory of finite groups.

Steps 2 and 3 already appear in the original definition of the 
Jones polynomial \cite{J2}, where each link is represented
as a closure of a braid, the Jones polynomial is defined from
such a braid through certain representation theory, and then
it is proved that this polynomial is independent of a choice
of a braid for a fixed link.

This strategy should work, in principle, in any dimension, but
so far, most of the interesting constructions arising from
operator algebras are for dimension 3,
so we concentrate in this case in this survey.

There have been many constructions of such topological invariants
for 3-dimensional manifolds and two of them are particularly related
to operator algebras.  One is a construction of Turaev-Viro \cite{TV}
in a generalized form due to Ocneanu, and the other is the one by
Reshetikhin-Turaev \cite{RT}.  For these two, the triple of
operator algebraic systems, representation theoretic data, and
the topological invariants in each case is listed as in Table \ref{tab1}.

\begin{table}
\centering
\caption{Topological invariants arising from operator algebras}
\label{tab1}
\begin{tabular}{lll}
\hline\noalign{\smallskip}
Operator Algebras & Representation Theory
& Combinatorial Construction  \\
\noalign{\smallskip}\hline\noalign{\smallskip}
Subfactors & Quantum $6j$-symbols & TVO invariants \\
\noalign{\smallskip}\hline\noalign{\smallskip}
Nets of factors on $S^1$ & Braided tensor categories & 
RT invariants \\
\noalign{\smallskip}\hline
\end{tabular}
\end{table}

Since both operator algebras and (topological) quantum field
theory are of infinite dimensional nature, one expects a
direct and purely infinite dimensional construction of the
latter from the former, but such a construction has not
been known yet.  All the constructions below go through representation
theoretic combinatorial data who ``live in'' finite dimensional
spaces, so one could eliminate the initial infinite dimensionality
entirely, if one is interested in only new constructions and
computations of topological invariants of 3-dimensional manifolds.
Still, the infinite dimensional framework
of operator algebras is useful, as we see below, even in such
a case, because it gives a conceptually convenient working place
for various constructions and computations.

We also mention one reason we operator algebraists are interested
in this type of theory, even purely from a viewpoint of operator algebras.
Classification theory is a central topic in theory of operator
algebras, and representation theory gives a very important invariant
for classification.  Since a series of great works of A. Connes in
1970's, it is believed that under some nice analytic condition,
generally called ``amenability'', a certain representation theory
should give a complete invariant of operator algebraic systems,
such as operator algebras themselves, group actions on them, 
or certain families of operator algebras.  For this reason,
studies of representation theories in operator algebraic theory
are quite important since old days of theory of operator algebras.
What is new after the emergence of the Jones theory is that
the representation theory now has a ``quantum'' nature, whatever
it means.

The author thanks R. Longo, N. Sato, and H. Wenzl for comments on this
manuscript.

\section{Turaev-Viro-Ocneanu invariants}

Here we review the Turaev-Viro-Ocneanu invariants of 3-dimensional
manifolds.  The book \cite{EK} is a basic reference.

Our operator algebra here is a so-called {\sl von Neumann algebra},
which is an algebra of bounded linear operators on a certain
Hilbert space that is closed under the $*$-operation and the
strong operator topology.  (Here we consider only infinite
dimensional separable Hilbert spaces, though a general theory
exists for other Hilbert spaces.)  Requiring closedness under the
weak operator topology, we obtain the same class of operator algebras.
If we use a norm topology, we have a wider class of operator
algebras called {\sl $C^*$-algebras}.  Although von Neumann algebras
give a subclass of $C^*$-algebras, it is not very useful,
except for some elementary aspects of the theory, to regard
a von Neumann algebra as a $C^*$-algebra, because a von Neumann
algebra is far from being a ``typical'' $C^*$-algebra.  For example,
most of natural $C^*$-algebras are separable, as Banach spaces, but
von Neumann algebras are never separable, unless they are finite
dimensional.  We assume, as usual, that a von Neumann algebra
contains the identity operator, which is the unit of the algebra.
A commutative $C^*$-algebra having a unit is the algebra of all
the continuous functions on a compact Hausdorff space, and a
commutative von Neumann algebra is the algebra of $L^\infty$-functions
on a measure space.  This gives a reason for a basic idea that
a general $C^*$-algebra is a ``noncommutative topological space''
and a general von Neumann algebra is a ``noncommutative measure
space''.   A finite dimensional $C^*$- or von Neumann algebra is
isomorphic to a finite direct sum of full matrix algebras $M_n({\bf C})$.
In this paper, we deal with only {\sl simple} von Neumann algebras
in the sense that they have only trivial two-sided closed ideals
in the strong or weak operator topology.  This simplicity is
equivalent to triviality of the center of the algebra, and we
call such a von Neumann algebra a {\sl factor}, rather than a simple
von Neumann algebra.

In the Murray-von Neumann classification, factors are
classified into type I, type II$_1$, type II$_\infty$, and 
type III.  Factors of type I are simply all the bounded linear
operators on some Hilbert space, and they are not interesting
for the purpose of this survey.  We are interested in factors
of type II$_1$ in the following two sections and those of type III in
the last section.  Although technical details on these factors
are not necessary for conceptual understanding of the theory,
we give brief explanations on how to construct such factors.

We start with a countable group $G$.  The (left) regular
representation gives a unitary representation of $G$ on
the Hilbert space $\ell^2(G)$.  We consider the von Neumann
algebra generated by its image.  If the group $G$ is commutative,
the resulting von Neumann algebra is isomorphic to $L^\infty(\hat G)$.
If the group $G$ is ``reasonably noncommutative'' in an appropriate
sense, the resulting von Neumann algebra is a factor of type II$_1$.
One example of such a group is that of all permutations of a countable
set that fix all but finite elements.

Another construction of a factor arises from an infinite tensor
product of the $n\times n$-matrix algebra $M_n({\bf C})$.  We can
define such an infinite tensor product in an appropriate sense, 
and then this infinite dimensional algebra has a natural representation
on a separable Hilbert space.  The von Neumann algebra generated
by its image is a type II$_1$ factor and these are all isomorphic,
regardless $n$.  This infinite tensor product also has many other
representations on Hilbert spaces and ``most'' of them generate
factors of type III.

The most natural starting point of a representation theory for
factors is certainly a study of all representations of a fixed
factor on Hilbert spaces.  (A factor is an algebra of operators 
on a certain Hilber space by definition, but we consider
representations on other Hilbert spaces.  In our setting, it
is enough to consider only representations on infinite dimensional
separable Hilbert spaces.)  We certainly have a natural notion 
of unitary equivalence for representations of factors of type II$_1$
or III, but this notion is not particularly interesting, as follows.
Such  representations are never irreducible, and for a
fixed  type II$_1$ factor, we can classify representations completely,
up to unitary equivalence, with a single invariant, called a coupling
constant, due to Murray and von Neumann, having values in $(0,\infty]$.
(This invariant produces the Jones index as below, and produces
something deep in this sense, but the classification of representations
themselves is rather simple and classical.)
For factors of type III, the situations are even simpler; they are
all unitarily equivalent for a fixed type III factor.

A representation of a factor can be regarded as a (left) module
over a factor, trivially.  It was Connes who realized first that the
right setting for studying representation theory of factors is to
study {\sl bimodules}, rather than modules.  That is, we consider
two factors $M$ and $N$, which could be equal, and study
a Hilbert space $H$ which is a left $M$-module and a right $N$-module
with the two actions commuting.
We call such $H$ an $M$-$N$ bimodule and write ${}_M H_N$.
The situation where both $M$ and $N$ are of type II$_1$ is
technically simpler.  We have natural notions of irreducible
decomposition, dimensions having
values in $(0,\infty]$ which are defined in terms of the
coupling constants, contragredient bimodules,
and relative tensor products.  For example, for two bimodules
${}_M H_N$ and ${}_N K_P$, we can define an $M$-$P$ bimodule
${}_M H\otimes_N K_P$ and the dimension is multiplicative.
For a factor $M$, the algebra $M$ itself trivially
has the left and right actions of $M$, so it has a bimodule
structure, but this $M$ is not a Hilbert space.  We have
a natural method to put an inner product on $M$ and complete it,
and in this way, we obtain an $M$-$M$ bimodule.  By an abuse of
notation, we often write ${}_M M_M$ for this bimodule, by
ignoring the completion.  This bimodule
has dimension one, and plays a role of a trivial representation.
In this way, our representation theory is quite analogous to that
of a compact group.
Connes used a terminology {\sl correspondences} rather than
bimodules.  See \cite{P1} for a general theory on bimodules.

Jones initiated studies of inclusions of factors $N\subset M$
in \cite{J1}.  Such $N$ is called a {\sl subfactor} of $M$.
By an abuse of
terminology, the inclusion $N\subset M$ is often called a
subfactor.  Technically simpler situations are that both
$M$ and $N$ are of type II$_1$.  Then we have an $M$-$M$
bimodule ${}_M M_M$ as above, and we restrict the left
action to the subalgebra $N$ to obtain ${}_N M_M$.  The
dimension of this bimodule is called the {\sl Jones
index} of the subfactor $N\subset M$ and written as
$[M:N]$.  (This terminology and notation come from an
analogy to a notion of an index of a subgroup.)  Jones
proved in \cite{J1} an astonishing statement that
this index takes values in $\{4\cos^2(\pi/n)\mid
n=3,4,5\dots\}\cup[4,\infty]$ and all the values in
this set are realized.  This is in a sharp contrast
to the fact that the coupling constant of a type II$_1$
factor $M$ can take all values in $(0,\infty]$.
Jones introduced the {\sl basic construction} whose
successive uses produce an increasing sequence
$N\subset M\subset M_1\subset M_2\subset\cdots$
and using this, he introduced the {\sl higher relative
commutants} and the {\sl principal graph} for subfactors.
Although we do not give their definitions here, we
only mention that if the subfactor has index less than 4,
then the principal graph is one of the $A$-$D$-$E$ Dynkin
diagrams, as noted by Jones.  (See \cite[Chapter 9]{EK}
for precise definitions.)

It was Ocneanu \cite{O1} who realized that these invariants
and further finer structures related to them can be
captured by theory of bimodules and that they can
be characterized by a set of combinatorial axioms.  We explain 
his theory here.  See \cite[Chapter 9]{EK} for more details.
We start with a type II$_1$ subfactor $N\subset M$
with finite Jones index.  (If we have a finite index
and one of $N$ and $M$ is of type II$_1$, then the other
is also of type II$_1$ automatically.)
Ocneanu's idea was to develop a {\sl representation theory
for a pair $N\subset M$}.  We start with ${}_N M_M$ and
this plays a role of the fundamental representation.
We also have ${}_M M_N$ and make relative tensor products
such as ${}_N M \otimes_M M \otimes_N M_M$.  They are
not irreducible in general, so we make irreducible
decompositions.  We look at
all unitary equivalence classes of $N$-$N$ bimodules
arising in this way.  In general, we expect to have
infinitely many equivalence classes, but it sometimes
happens that we have only finitely many equivalence
classes.  This is the situation we are interested in,
and in such a case, we say that the subfactor $N\subset M$
has a finite depth.  (The terminology ``depth'' comes
from the way of Jones to write higher relative commutants.)
This finite depth condition is similar to rationality
condition in conformal field theory and quantum group
theory.  If we have a finite depth, we also have only
finitely many equivalence classes of irreducible $M$-$M$
bimodules arising in the above way. 
Note that a compact group has only finitely
many equivalence classes of irreducible unitary
representations if and only if the group is finite.
We assume the finite depth condition and fix a finite set
of representatives of equivalence classes of irreducible
$N$-$N$ bimodules arising as above from $N\subset M$.
Note that it contains a trivial bimodule, that for
each bimodule in the set, its contragredient bimodule
is equivalent to one in the set,  and that a relative
tensor product of two in the set decomposes into a finite
direct sum of irreducible bimodules each of which is
equivalent to one in the set.  We say such a finite
set of bimodules is a finite system of bimodules.
Choose three, not necessarily distinct,
irreducible $N$-$N$ bimodules $A, B, C$ in the system.
Then we can decompose $A\otimes_N B \otimes_N C$ in two
ways.  That is, we first decompose $A \otimes_N B$ in
one, and we first decompose $B\otimes_N C$ in the other.
In this way, we obtain the ``quantum'' version of the
classical $6j$-symbols which produce a complex number
from six bimodules and four intertwiners.  Such
quantum $6j$-symbols were known in the quantum group
theory, and Ocneanu found that a general system of bimodules
produce similar $6j$-symbols and that classical properties
such as the Frobenius reciprocity also holds in this setting.
Associativity of the relative tensor product gives a
so-called pentagonal relation as in the classical
setting.  This finite system of bimodules and quantum
$6j$-symbols are the combinatorial data arising from a
representation theory of a subfactor $N\subset M$.

Turaev and Viro \cite{TV} constructed topological invariants of
3-dimensional manifolds using the quantum $6j$-symbols
for the quantum group $U_q(sl_2)$ at roots of unity, and
Ocneanu realized that a generalized version of
this construction works for general
quantum $6j$-symbols arising from a subfactor of  finite Jones
index and finite depth
as above.  The construction goes as follows for a fixed
finite system of bimodules.  (See \cite[Chapter 12]{EK} for
more details.)

We first make a triangulation of a manifold.  That is,
we regard a manifold made of gluing
faces of finitely many tetrahedra so that we have
an empty boundary and compatible orientation.  Then we label
each of the six edges with bimodules in the system and each of
the four faces, triangles, with (co-)isometric intertwiners.  
When all the tetrahedra
are labeled in this way, the quantum $6j$-symbol produce a
complex number for each labeled tetrahedron.  This number
is simply a composition of the four intertwiners, up to
normalization arising from dimensions of the four bimodules.
(The composed intertwiners give a complex number because of
irreducibility of the bimodules.)  The well-definedness of
this number comes from the so-called tetrahedral symmetry
of quantum $6j$-symbols.  Then we multiply all these numbers 
over all the tetrahedra in the triangulation, and add these
products over all isometric intertwiners in an
orthonormal basis for each face and over all labeling of
edges with bimodules.  With an appropriate normalization
arising from dimensions of the bimodules, the resulting number
is a topological invariant of the original 3-dimensional
manifold.  In order to prove this topological invariance, 
one has to prove that the complex number is independent of
triangulations of a manifold.  The relations of two triangulation
of a manifold have been known by Alexander.  That is, one
triangulation is obtained from the other by successive
applications of finitely many local changes of triangulations,
called Alexander moves.  (This result of Alexander holds in
any dimension.)  Pachner has proved that a different set of
local moves also gives a similar theorem, and this set is more
convenient for our purpose.  That is, it is enough for us
to prove that the above complex number is invariant under
each of the Pachner moves.  This invariance follows from
properties of the quantum $6j$-symbols, such as the pentagon
relation.  So we conclude that the above complex number gives
a well-defined topological invariant of 3-dimensional
closed oriented manifolds.  If we reverse the orientation,
the topological invariant becomes the complex conjugate of 
the original value.  In the original setting of Turaev-Viro
\cite{TV} based on the quantum $6j$-symbols of $U_q(sl_2)$, the
resulting invariants are real, so they do not detect
orientations, but there is an example of a subfactor which
produces a non-real invariant for some manifold and thus
can detect orientations.  (Actually, the original construction
of Turaev-Viro \cite{TV} works without orientability.)  Also note
that in our setting, each intertwiner space has a Hilbert
space structure and each dimension of a bimodule, which 
is sometimes called a {\sl quantum dimension}, is positive.
Such a feature is called {\sl unitarity} of quantum $6j$-symbols,
and this unitarity does not necessarily hold in a purely algebraic
setting of quantum $6j$-symbols for quantum groups.
We can apply the same construction by using the system
of the $M$-$M$ bimodules instead of that of the $N$-$N$
bimodules, but the resulting invariant is the same.

A large class of subfactors are constructed with
methods related to classical theory of groups and Hopf
algebras, and their ``quantum'' counterparts, that is,
quantum group theory and
conformal field theory such as the Wess-Zumino-Witten models.
For such subfactors, we have various interesting studies
from an operator algebraic viewpoint, but if we are interested
only in resulting topological invariants through the above
machinery, they do not produce really new invariants.
It is, however, expected that we have much wider varieties
of subfactors in general.  One ``evidence'' for such
expectation is study of Haagerup \cite{Ha}.  By purely
combinatorial arguments, he found a list of candidates of
subfactors of finite depth in the index range $(4,3+\sqrt2)$,
and it seems that most of these are indeed realized.  None
of them seem to be related to conformal field theory
or today's theory of quantum groups so far.  Haagerup himself
proved that the first one in the list is indeed realized,
and Asaeda and he further proved that another in the list
is also realized in \cite{AH}.  The nature of topological
invariants arising from these two subfactors is not understood
yet, but we expect that they contain some interesting information.
Since the list of Haagerup is only for a small range of the
index values, we expect that we would have by far more 
examples of ``exotic'' subfactors as mentioned above, but an explicit
construction of even a single example is highly difficult.
We know almost nothing about topological meaning of invariants
arising from such subfactors.  Izumi \cite{I2} has some more
examples of such interesting subfactors.

\section{Reshetikhin-Turaev invariants}

Another construction of topological invariants due to
Reshetikhin-Turaev \cite{TV} requires a ``higher symmetry''
for combinatorial data arising from a representation theory.
This higher symmetry is called a {\sl modularity} of a tensor
category.  It is also called a {\sl nondegenerate braiding}.

Wenzl has a series of work \cite{W,TW,TW2,W2},
partly with V. G. Turaev, on related constructions, but here
we concentrate on two methods producing a modular tensor category from
a general operator algebraic representation theory.  One is within
subfactor theory, due to Ocneanu, and presented in this
section, and the other is due to
Longo, M\"uger and the author \cite{KLM}, explained in the
next section.

We first give a brief explanation on braiding.  In a
representation theory of a group, two tensor products
$\pi\otimes\sigma$ and $\sigma\otimes \pi$ are
obviously unitarily equivalent for two representations
$\pi$ and $\sigma$, but for two $N$-$N$ bimodules
$A$, $B$, we have no reason to expect that
$A\otimes_N B$ and $B\otimes_N A$ are equivalent, and
they are indeed not equivalent in general.  It is, however,
possible that for all $A$ and $B$ in a finite system,
we have equivalence of $A\otimes_N B$ and $B\otimes_N A$.
If we can choose isomorphisms of these two bimodules
in a certain compatible way simultaneously for all
bimodules in the system, we say that the system has
a {\sl braiding}.  See \cite{R} for more details, where
an equivalent, but slightly different formulation using
endomorphisms, rather than bimodules, is presented.

The isomorphism between $A\otimes_N B$ and $B\otimes_N A$
can be graphically represented as an overcrossing of
two wires labeled with $A$ and $B$, respectively.  Then
the assumption on ``compatibility'' implies, for example,
the Yang-Baxter equation, which represents the Reidemeister
move of type III as in Fig. \ref{R3}, where each crossing
represents an isomorphism and each hand side is a composition
of three such isomorphisms.

\thicklines
\unitlength 1mm
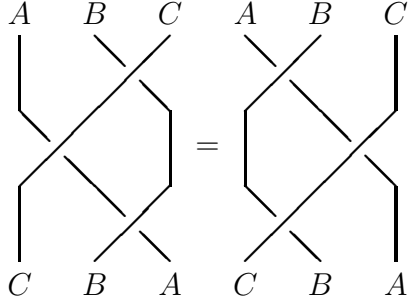
\begin{figure}[tb]
\begin{center}
\begin{picture}(70,50)
\put(10,10){\line(0,1){10}}
\put(10,20){\line(1,1){20}}
\put(20,10){\line(1,1){10}}
\put(30,20){\line(0,1){10}}
\put(30,30){\line(-1,1){4}}
\put(20,40){\line(1,-1){4}}
\put(30,10){\line(-1,1){4}}
\put(20,20){\line(1,-1){4}}
\put(20,20){\line(-1,1){4}}
\put(10,30){\line(1,-1){4}}
\put(10,30){\line(0,1){10}}
\put(40,10){\line(1,1){20}}
\put(60,30){\line(0,1){10}}
\put(50,10){\line(-1,1){4}}
\put(40,20){\line(1,-1){4}}
\put(40,20){\line(0,1){10}}
\put(40,30){\line(1,1){10}}
\put(60,10){\line(0,1){10}}
\put(60,20){\line(-1,1){4}}
\put(50,30){\line(1,-1){4}}
\put(50,30){\line(-1,1){4}}
\put(40,40){\line(1,-1){4}}
\put(10,7){\makebox(0,0){$C$}}
\put(20,7){\makebox(0,0){$B$}}
\put(30,7){\makebox(0,0){$A$}}
\put(40,7){\makebox(0,0){$C$}}
\put(50,7){\makebox(0,0){$B$}}
\put(60,7){\makebox(0,0){$A$}}
\put(10,43){\makebox(0,0){$A$}}
\put(20,43){\makebox(0,0){$B$}}
\put(30,43){\makebox(0,0){$C$}}
\put(40,43){\makebox(0,0){$A$}}
\put(50,43){\makebox(0,0){$B$}}
\put(60,43){\makebox(0,0){$C$}}
\put(35,25){\makebox(0,0){$=$}}
\end{picture}
\end{center}
\caption{Yang-Baxter equation}
\label{R3}
\end{figure}

In representation theory of groups, the tensor product
operation is trivially commutative in the above sense.
This is ``too commutative'' in the sense that we have
no distinction between an overcrossing and an undercrossing
in the above graphical representation, and this is not very
useful for construction of topological invariants,
obviously.  So, in
order to obtain an interesting topological invariant, an
overcrossing and an undercrossing must be ``sufficiently
different''.  Such a condition is called nondegeneracy of
the braiding.  This condition can be also formulated in the
language of tensor categories, and then it is called a
modularity of the tensor category.  A non-degenerate
braiding, or a modular tensor category, produces a unitary
representation of the modular group $SL(2,{\bf Z})$.
 
We first explain how to obtain such a nondegenerate braiding
in subfactor theory.  We start with
a subfactor $N\subset M$ with finite Jones index and
finite depth.  Then Ocneanu has found a construction of
a new subfactor from this subfactor, which is called
the asymptotic inclusion \cite{O1}.  He realized that
the system of bimodules for this new subfactor has a nondegenerate
braiding and it can be regarded as the ``quantum double'' of
the original system of $N$-$N$ (or $M$-$M$) bimodules arising
from the subfactor $N\subset M$.  Note that the original
system of $N$-$N$ bimodules and that of $M$-$M$ bimodules are
not isomorphic in general, but they have the same ``quantum
double'' system of bimodules.  Popa has a more general construction
of this type, called the symmetric enveloping algebra \cite{P2}.
Longo-Rehren \cite{LR} has found the essentially same
construction as the asymptotic inclusion in the setting of
algebraic quantum field theory.  See \cite[Chapter 12]{EK} for
more details on the asymptotic inclusion and \cite{I1,I2} for
detailed analysis based on the Longo-Rehren approach.

Suppose we have a nondegenerate braiding.  It is also known that
such a braiding can arise from quantum groups or conformal field
theory.  Reshetikhin-Turaev \cite{RT} has constructed a topological
invariant of 3-dimensional manifold from such a system.
First we draw a picture of a link on a plane.  This has various
overcrossings and undercrossings.  We label each connected component
with an irreducible bimodule in the system, then each crossing
gives an isomorphism arising from the braiding.  Then this labeled
picture produces a complex number as a composition of these
isomorphisms.  This is an invariant of  ``colored links'', where
coloring means labeling of each component with an irreducible
bimodule.  Actually, this number is not invariant under the
Reidemeister move of type I, and it is invariant under only the
Reidemeister moves of type II and type III, so this is not a
topological invariant of colored links, but it gives a
``regular isotopy'' invariant of colored links, for which
invariance under the Reidemeister moves of type II and type III
is sufficient.  Then we sum these complex numbers over all
possible colorings, with appropriate normalizing weights arising
from dimensions of the bimodules.  In this way, we obtain a
complex number from a planar picture of a link.  There is a
method to construct a 3-dimensional oriented closed
manifold from such a planar picture of a link, called the
{\sl Dehn surgery}.  Roughly speaking, we embed a link in the
3-sphere, and remove a tubular neighbourhood, consisting of
a disjoint union of solid tori, from the 3-sphere, and then
put back the solid tori in a different way.
Different links can produce the same
3-dimensional manifolds, but again, it is known that
in such a case, the two links can be transformed from one
to the other with successive applications of local moves.
Such moves are called Kirby moves.  Reshetikhin and Turaev
have proved that nondegeneracy of the braiding implies
invariance of the above complex number, the weighted
sum of colored link invariants, under Kirby moves, thus
we obtain a topological invariant of 3-dimensional 
manifolds in this way.  Reshetikhin and Turaev considered
an example arising from the quantum groups $U_q(sl_2)$ at
roots of unity, but the general machinery applies to any
nondegenerate braiding.  See the book \cite{T} for more
details on this construction.

So, starting with a subfactor with finite Jones index and
finite depth, we have two topological invariants of
3-dimensional manifolds.  One is the Turaev-Viro-Ocneanu
invariant arising from the system of $N$-$N$ bimodules.
The other is the Reshetikhin-Turaev invariant of the
``quantum double'' system of the original system of
$N$-$N$ bimodules.  It is quite natural to investigate
the relation between these two invariants.  Sato, Wakui
and the author proved in \cite{KSW} that these two invariants
coincide.  Ocneanu \cite{O3} has also announced such
coincidence and it seems to us that his method is different
from ours. Sato and Wakui \cite{SW} also made explicit
computations of this invariant for various concrete examples
of subfactors and manifolds, based on Izumi's explicit
computations of the representations of the modular group
arising from some subfactors, including the ``exotic''
one due to Haagerup, in \cite{I2}.
 
Another computation of topological invariants arising from
subfactors is based on $\alpha$-induction
\cite{LR,X1,BE,BEK1,BEK2}.  This method, in particular, 
produces subfactors with principal graphs $D_{2n}$,
$E_6$, and $E_8$, and the corresponding Turaev-Viro-Ocneanu
invariants can be computed once we have a description of the
``quantum doubles'' by \cite{KSW}, and these quantum doubles
were computed in \cite{BEK3}.  (Also see \cite{O3}.)  
This $\alpha$-induction
is also related to theory of modular invariants \cite{CIZ}.
See \cite{BE,BEK1,BEK2,KL,KL2} for more on this topic.

\section{Algebraic quantum field theory}

Another occurrence of nondegenerate braiding in theory of
operator algebras is in algebraic quantum field theory \cite{H},
which has its own long history.
This theory is an approach to quantum field theory based on
operator algebras.  That is, in each bounded region on a
spacetime, we assign a von Neumann algebra on a fixed Hilbert
space.  We think that each such von Neumann algebra is generated
by observable physical quantities in the bounded region in
the spacetime.  In this way, we think that this family of
von Neumann algebras parametrized by bounded regions gives
a mathematical description of a physical theory.  We often
restrict bounded regions to those of a special form.  We impose
a physically natural set of axioms on this family of von Neumann
algebras and make a mathematical study of such axiomatized
systems.  A spacetime of any dimension is allowed in this
axiomatized approach, and the four dimensional case was studied
originally for an obvious physical reason.  These studies of
Doplicher-Haag-Roberts \cite{DHR} and Doplicher-Roberts 
\cite{DR1,DR2} have been quite successful.  Recently, it has
been realized that this theory in lower dimensional spacetime
has quite interesting mathematical structures.  Two-dimensional
case has caught much attention in connection to conformal
field theory and one-dimensional case also naturally appears 
in a ``chiral'' decomposition of a two-dimensional theory.
Mathematical structures of one-dimensional theory was studied
in \cite{FRS}.  In one-dimensional case, our ``spacetime''
is simply ${\bf R}$ and a bounded region is simply a bounded
interval.  It is often convenient to compactify the
space ${\bf R}$ to obtain $S^1$ and consider ``intervals''
contained in $S^1$.  In this setting, our mathematical
structure is a family of von Neumann algebras on a fixed
Hilbert space parameterized by intervals in $S^1$.  We impose
a set of axioms.  For example, one axiom requires that we
have a larger von Neumann algebra for a larger interval.  Another
axiom requires ``covariance'' of the theory with respect to a
projective unitary representation of a certain group of
the ``spacetime symmetry''.  We also have an axiom on
``locality'' which says if two regions are spacelike
separated, then the corresponding von Neumann algebras
mutually commute.  Another requires existence of
a ``vacuum'' vector in the Hilbert space.  Positivity of
energy in the sense that a certain self-adjoint operator
is positive is also assumed.  See \cite{GL,KLM}
for a precise description of the set of axioms.  (Actually
the main results in \cite{KLM} hold under  a weaker set of
axioms, but we do not go into details here.)  Under the
usual set of axioms, each von Neumann algebra for an interval
becomes a factor of type III, so we call such a family a 
{\sl net of factors}.  Now the index set of intervals on
the circle $S^1$ is not directed with respect to inclusions,
since the entire circle is not allowed as an interval, so
it is not appropriate to call such a family a {\sl net}, but
this terminology has been commonly used.

This family is our operator algebraic system and we consider
a representation of such a family of von Neumann algebras.
Such an idea is due to Doplicher-Haag-Roberts \cite{DHR}
and is called the DHR theory.  We have a natural notion of
irreducibility, dimensions, and tensor products for such
representations.  Note that we do not have an obvious
definition of tensor products for two representations of
such a net of factors.  The key idea was that the tensor
product operation is realized through compositions of
endomorphisms.  Also the dimension in the usual
sense is always infinite.  So it was highly nontrivial to
obtain sensible definitions of the tensor product and
the dimension.  This work is much
older than the subfactor theory in the previous section, 
and its similarity to subfactor theory was soon
recognized in \cite{L} in a precise form.

In this way, we have a representation theory for a net
of factors.  A tensor product operation is ``too commutative''
for higher dimensional spacetime, but
in dimensions one and two, it has an appropriate level
of commutativity, and naturally produces a braiding.
(See \cite{FRS} for example.)  So we have two problems for
getting a modular tensor category from such a 
representation of a net of factors on $S^1$.  One is
whether we have only finitely many equivalence classes
of irreducible representations or not.  The other is
whether the braiding is nondegenerate or not.  In
\cite{KLM}, Longo, M\"uger and the author have found a
nice operator algebraic condition that implies positive
answers to these two problems and we introduced the
terminology ``complete rationality'' for this notion.
One of the key conditions for this notion is finiteness
of a certain Jones index.  Note that in subfactor
theory in the previous section, our ``family of operator
algebras'' has only two factors $N$ and $M$, and its
representation theory produced a tensor category,
without a braiding in general.  Now our ``family of operator
algebras'' is a net of factors and has continuously many factors
with more structures, and its representation theory produces
a braided tensor category.

Xu has proved in \cite{X2} that the $SU(N)_k$-nets
corresponding to the WZW-models $SU(N)_k$ are completely
rational.  Xu worked on coset models in the setting
of nets of factors on $S^1$ in \cite{X3}, and obtained
several interesting examples.  He then studied in
\cite{X4} about topological invariants arising from
these nets, which seems to be quite interesting
topologically.  He also worked on orbifold models
in this context in \cite{X5}.
Finally, we also note that complete rationality is also important
in classification theory of nets of factors as in
\cite{KL,KL2}.


\begin{thebibliography}{99}

\bibitem{AH}
M. Asaeda, U. Haagerup:
Commun. Math. Phys. \textbf{202}, 1--63 (1999)

\bibitem{A}
M. F. Atiyah: Publ. Math. I.H.E.S.  \textbf{68}, 175--186 (1989)

\bibitem{BE}
J. B\"ockenhauer, D. E. Evans:
Commun. Math. Phys. \textbf{197}, 361--386 (1998)
II \textbf{200},  57--103 (1999)
III \textbf{205}, 183--228 (1999)

\bibitem{BEK1}
J. B\"ockenhauer, D. E. Evans, Y. Kawahigashi:
Commun. Math. Phys.  \textbf{208}, 429--487 (1999)

\bibitem{BEK2}
J. B\"ockenhauer, D. E. Evans, Y. Kawahigashi:
Commun. Math. Phys. \textbf{210}, 733--784 (2000)

\bibitem{BEK3}
J. B\"ockenhauer, D. E. Evans, Y. Kawahigashi:
Publ. RIMS, Kyoto Univ. \textbf{37}, 1--35 (2001)

\bibitem{CIZ}
A. Cappelli, C. Itzykson, J.-B. Zuber:
Commun. Math. Phys.  \textbf{113}, 1--26 (1987)

\bibitem{DHR}
S. Doplicher, R. Haag, J. E. Roberts:
I Commun. Math. Phys.  \textbf{23}, 199--230 (1971)
II \textbf{35}, 49--85 (1974)
 
\bibitem{DR1}
S. Doplicher, J. E. Roberts:
Ann. Math. \textbf{130}, 75--119 (1989)

\bibitem{DR2}
S. Doplicher, J. E. Roberts:
Invent. Math. \textbf{98}, 157--218 (1989)

\bibitem{EK}  D. E. Evans, Y. Kawahigashi:
\textit{Quantum symmetries on operator algebras},
(Oxford University Press, Ofxord, 1998)

\bibitem{FRS}
K. Fredenhagen, K.-H. Rehren, B. Schroer:
I. Commun. Math. Phys.  \textbf{125}, 201--226 (1989)
II. Rev. Math. Phys. \textbf{Special issue}, 113--157 (1992)

\bibitem{GL}
D. Guido, R. Longo:
Commun. Math.  Phys. \textbf{181}, 11--35 (1996)

\bibitem{H} R. Haag:
\textit{Local Quantum Physics},
(Springer-Verlag, Berlin-Heidelberg-New York, 1996)

\bibitem{Ha} U. Haagerup:
Principal graphs of subfactors in the index range
$4< 3+\sqrt2$. In:
\textit{Subfactors}, ed by H. Araki, et al.
(World Scientific, 1994) pp 1--38

\bibitem{I1}
M. Izumi: Commun. Math. Phys.  \textbf{213}, 127--179 (2000)

\bibitem{I2}
M. Izumi: Rev. Math. Phys. \textbf{13}, 603--674 (2001)

\bibitem{J1}
V. F. R. Jones: Invent. Math. \textbf{72}, 1--25 (1983)

\bibitem{J2}
V. F. R. Jones: Bull. Amer. Math. Soc. \textbf{12}, 103--112 (1985)

\bibitem{KL}
Y. Kawahigashi, R. Longo:
math-ph/0201015, to appear in Ann. Math.

\bibitem{KL2}
Y. Kawahigashi, R. Longo:
math-ph/0304022, to appear in Commun. Math. Phys.

\bibitem{KLM}
Y. Kawahigashi, R. Longo, M. M\"uger:
Commun. Math. Phys.  \textbf{219}, 631--669 (2001)

\bibitem{KSW}
Y. Kawahigashi, N. Sato, M. Wakui: math.OA/0208238

\bibitem{L}
R. Longo: 
I Commun. Math. Phys. \textbf{126}, 217--247 (1989)
II Commun. Math. Phys. \textbf{130}, 285--309 (1990)

\bibitem{LR}
R. Longo, K.-H. Rehren: Rev. Math. Phys.
\textbf{7}, 567--597 (1995)

\bibitem{M}
M. M\"uger: math.CT/0111205

\bibitem{O1}
A. Ocneanu:
Quantized group, string algebras and Galois theory for algebras.
In \textit{Operator algebras and applications, Vol. 2},
ed D. E.  Evans and M. Takesaki, 
(Cambridge University Press, Cambridge, 1988) pp 119--172

\bibitem{O2}
A. Ocneanu: Chirality for operator algebras. In:
\textit{Subfactors}, ed by H. Araki, et al.
(World Scientific, 1994) pp 39--63

\bibitem{O3}
A. Ocneanu: Operator algebras, topology and subgroups of quantum
symmetry -- construction of subgroups of quantum groups --
(written by S. Goto and N. Sato).  In: 
\textit{Taniguchi Conference in Mathematics Nara '98}
Adv. Stud. Pure Math. \textbf{31},
(Math. Soc. Japan, 2000) pp 235--263

\bibitem{P1} S. Popa: Correspondences, preprint 1986

\bibitem{P2}
S. Popa: Math. Res. Lett.
\textbf{1}, 409--425 (1994)

\bibitem{R}
K.-H. Rehren:
Braid group statistics and their superselection rules.
In: \textit{The algebraic theory of superselection sectors,
Palermo, 1989}, World Scientific Publishing (1990) pp 333--355

\bibitem{RT}
N. Reshetikhin, V. G. Turaev:
Invent. Math. \textbf{103}, 547--597 (1991)

\bibitem{SW}
N. Sato and M. Wakui:
math.OA/0208242, to appear in J. Knot Theory Ramif.

\bibitem{T}
V. G. Turaev, \textit{Quantum Invariants of Knots and $3$-manifolds}, 
(Walter de Gruyter, 1994)

\bibitem{TV}
V. G. Turaev, O. Ya Viro: Topology \textbf{31}, 865--902 (1992)

\bibitem{TW}
V. G. Turaev, H. Wenzl:
Internat. J. Math. \textbf{4}, 323--358 (1993)

\bibitem{TW2}
V. G. Turaev, H. Wenzl:
Math. Ann. \textbf{309}, 411--461 (1997)

\bibitem{W}
H. Wenzl:
Invent. Math. \textbf{114}, 235--275 (1993)

\bibitem{W2}
H. Wenzl:
J. Amer. Math. Soc. \textbf{11}, 261--282 (1998)

\bibitem{X1}
F. Xu:
Commun. Math. Phys. \textbf{192}, 347--403 (1998)

\bibitem{X2}
F. Xu:
Commun. Contemp. Math. \textbf{2}, 307--347 (2000)

\bibitem{X3}
F. Xu:
Commun. Math. Phys. \textbf{211},  1--44 (2000)

\bibitem{X4}
F. Xu:
math.GT/9907077

\bibitem{X5}
F. Xu:
Proc. Nat. Acad. Sci. U.S.A.
\textbf{97}, 14069--14073 (2000)

\end{thebibliography}
\end{document}